# ON THE WEAK KRULL SYMMETRY OF A NOETHERIAN RING.

C.L.Wangneo

Jammu,J&K,India,180002

(E-mail:-wangneo.chaman@gmail.com )

**Abstract:-** We define when a noetherian ring R is called a right ( or a left) weakly krull symmetric ring . We then prove that if R is a right ( or a left ) krull homogenous ring then R is a right ( or a left ) weakly krull symmetric ring . This result modifies the main result of [2] . The key terms introduced in this paper are of independent interest .

## Introduction :-

In this paper for a noetherian ring R with nilradical N(R) we define when R is called a right ( or a left) weakly krull symmetric ring . We then state results that culminate in our main theorem which states that if R is a right ( or a left) krull homogenous ring then R is a right ( or a left ) weakly krull symmetric ring . This result modifies the main result of [2] . The key terms defined in this paper are that for a noetherian ring R with nilradical N(R) we introduce definitions of ideals A(R) and B(R) of R and certain subsets $\Lambda(R)$ and $\Lambda'(R)$ of min. spec.R . We also define semiprime ideals P(R) and Q(R) of R which we call as the right and the left krull homogenous parts respectively of the nilradical N(R) of R . The above key terms introduced are of independent interest . The paper is divided into two sections . In section (1) we introduce the definitions of ideals A(R) and B(R) of a noetherian ring R ,

where $A(R)$ is the sum of all right Ideal $I$ of $R$ such that $|I|_r < |R|_r$, and $B(R)$ is the sum of all left ideal $K$ of $R$ with $|K|_l < |R|_l$. We then prove some basic results regarding the ideals $A(R)$ and $B(R)$.

Next in section (2) we define for a noetherian ring $R$ certain subsets $\Lambda(R)$ and $\Lambda'(R)$ of min. spec.$R$ as follows;

$\Lambda(R) = \{ P_i \text{ in spec. } (R) \, / \, |R/P_i|_r = |R|_r \}$ and

$\Lambda'(R) = \{ Q_i \text{ in spec. } (R) \, / \, |R/Q_i|_l = |R|_l \}$.

We then use these definitions and the results of sections (1) and (2) to prove our main theorem stated above which is a modification of the main theorem of [2] and which essentially states that if $R$ is a right (or a left) krull homogenous ring then $R$ is a right (or a left) weakly krull symmetric ring.

**Notation and Terminology:-**

Throughout in this paper a ring is meant to be an associative ring with identity which is not necessarily a commutative ring. We will throughout adhere to the same notation and definitions as in [1]. Thus by a noetherian ring $R$ we mean that $R$ is both a left as well as a right noetherian ring. By a module $M$ over a ring $R$ we mean that $M$ is a right $R$-module unless stated otherwise. For the basic definitions regarding noetherian modules and noetherian rings and all those regarding Krull dimension, we refer the reader to [1]. Recall from [1], for a ring $R$ a right $R$-module $M$ is said to be right (or left) Krull-homogenous if every non zero right (or left) $R$-sub-module $N$ of $M$ has same right (or left) Krull dimension as that of $M$. Moreover we will use the following terminology throughout. If $R$ is a ring then we denote by Spec.$R$, the set of prime ideals of $R$ and by min.Spec.$R$, the set of minimal prime ideals of R . Again if R is a ring and M is a right R module We denote by r(T) the right annihilator of a subset T of M and by l(T) the left annihilator of a subset T of W in case W is a left R module . Also $|M|_r$ denotes the right Krull dimension of the right R-module M if it exists and $|W|_l$ denotes the left Krull dimension of the left R module W if it exists . For two subsets A and B of a given set A ≤ B means B contains A and A<B denotes A ≤ B but A≠B. Also for two sets A and B , A⊄B denotes the set B that does not contain the subset A. For an ideal A of R , c(A) denotes the set of elements of R that are regular modulo A .

**Section (1) :-** In this section we first introduce the definitions of the ideals A(R) and B(R) for a noetherian ring R . We then prove some basic results about them .

**Definition(1.1) :-** Let R be a noetherian ring. Denote by A(R) the sum of all right Ideal I of R such that $|I|_r < |R|_r$ . Similarly denote by B(R) the sum of all left ideal K of R with $|K|_l < |R|_l$ . If there is no confusion regarding the underlying ring R we may write A and B for A(R) and B(R) respectively.

We next state the obvious lemma (1.2) below;

**Lemma (1.2) :-** Let R be a Noetherian ring. Let J be an ideal of R which is faithful as a right ideal of R . Then $|J|_r = |R|_r$ .

**Proposition(1.3):-** Let R be a Noetherian ring. Let A and B be as in definition(1.1) above. Then the following hold true ;

(a) A is an ideal of R and is the unique largest right ideal of R with $|A|_r < |R|_r$ and $|R/A|_r = |R|_r$ .

(b)   R/A is a right k-homogenous ring. Also, if $A \neq 0$, then $r(A)$ is a nonzero proper ideal of R and $lr(A)=A$. Moreover $|A|_r = |R/r(A)|_r < |R|_r$.

**Proof:-** (a) Note first that R is right Noetherian ring implies that A is a finitely generated right ideal of R. Hence $|A|_r < |R|_r$. Also it is not difficult to see that A is an ideal of R and is the unique largest right ideal such that $|A|_r < |R|_r$. Moreover, it is obvious that $|R/A|_r = |R|_r$.

(b) It is easy to see that R/A is a right K-homogenous ring. We now show that $lr(A)=A$. So, let $A_1 = lr(A)$. By lemma (1.2) above it follows that $|A_1|_r = |R/r(A)|_r < |R|_r$. Hence, $A_1 \leq A$. Since already $A \leq A_1$, so $A = A_1$. Thus $lr(A)=A$. Again by lemma (1.2) above it follows that $|A|_r = |R/r(A)|_r < |R|_r$.

The left analogue of the above proposition(1.3) is a similar proposition regarding the ideal B and the ring R/B. It states the following;

**Proposition(1.4) :-** Let R be a Noetherian ring. Let A and B be as in definition(1.1) above. Then the following hold true ;

(a) B is an ideal of R and is the unique largest left ideal of R with $|B|_l < |R|_l$ and $|R/B|_l = |R|_l$.

(b) R/B is a left k-homogenous ring. Also, if $B \neq 0$, then $l(B)$ is a nonzero proper ideal of R and $rl(B)=B$. Moreover,

$|B|_l = |R/l(B)|_l < |R|_l$.

We now state Proposition (1.5) below that stengthens the conclusion of proposition (1.4) for a noetherian right krull homogenous ring.

**Proposition(1.5) :-** Let R be a Noetherian right k-homogenous ring . Let B be the largest left ideal of R with $|B|_l < |R|_l$ . Then the following hold true ;

(a) B is an ideal of R such that $|R/B|_l = |R|_l$ . Also R/B is a left k-homogeneous ring .

(b) The ring R/B is such that $|R/B|_r = |R|_r$ .

(c) Moreover R/B is a left as well as a right k-homogeneous ring with $|R/B|_r = |R|_r$ and $|R/B|_l = |R|_l$ .

**Proof:-** (a) is true by the left analogue of proposition (1.3) above .

(b) Observe that (b) is true if B = 0 . So, assume that B ≠ 0. We first show that $|R/B|_r = |R|_r$ . By the left analogue of the previous proposition , we get that $l(B) \neq 0$ and $l(B)$ is a proper ideal of R . Moreover, $|l(B)|_r \leq |R/B|_r$. Since R is right k-homogeneous thus we get that $|l(B)|_r = |R|_r$ and hence $|R/B|_r = |R|_r$.

(c) We already know from (a) above that R/B is a left krull homogenous ring . We now show that R/B is also a right k-homogeneous ring . Suppose not, then there exists a nonzero ideal $A_1/B$ of R/B such that $|A_1/B|_r < |R|_r$ and $A_1/B$ is largest such right ideal of R/B . In view of the previous proposition, we get that if Y/B = r-ann.$(A_1/B)$, then Y/B is nonzero ideal and $|R/Y|_r < |R|_r$ . Now, $A_1Y \leq B$ implies that $l(B) A_1 Y = 0$ . To show that R/B is also a right k-homogeneous ring , we first show that we can not have $l(B) A_1 \neq 0$ . To see this assume $l(B) A_1 \neq 0$ , then since $l(B)A_1Y = 0$, we must have that $|l(B) A_1|_r \leq |R/Y|_r$. But from above since $|R/Y|_r < |R|_r$ so we must have that $|l(B) A_1|_r < |R|_r$. This contradicts that R is a right k-homogenous ring. Next we show that we can not have $l(B) A_1 = 0$. For if $l(B) A_1 = 0$, then we get that $|A_1|_l \leq |R/l(B)|_l < |R|_l$. This implies that that $A_1 \leq B$ which is against our hypothesis since we have assumed that that $A_1/B$ is a non zero ideal. Hence, all in all we must have that $A_1 = B$ and so $R/B$ is a right k-homogeneous ring.

The statement of the left analogue of proposition(1.5) above is the following proposition;

**Proposition(1.6)** :- Let R be a Noetherian left k-homogenous ring. Let A be the largest right ideal of R with $|A|_r < |R|_r$. Then the following hold true;

(a) A is an ideal of R such that $|R/A|_r = |R|_r$. Also R/A is a right k-homogeneous ring.

(b) The ring R/A is such that $|R/A|_l = |R|_l$.

(c) Moreover R/A is a left as well as a right k-homogeneous ring with $|R/A|_l = |R|_l$ and $|R/A|_r = |R|_r$.

**Section (2)** :- In this section we introduce for a noetherian ring R, the definitions of certain subsets of min. spec.(R) denoted by $\Lambda(R)$ and $\Lambda'(R)$ respectively. Next we define when R is called a right ( or a left) weakly krull symmetric ring, a weakly krull symmetric ring and a weakly ideal krull symmetric ring. We then state and prove results that culminate in our main theorem.

**Definition(2.1) :-** (i) Let R be a Noethrian ring . Define two subsets of spec. R denoted by $\Lambda(R)$ and $\Lambda'(R)$ as follows ;

$\Lambda(R) = \{ P_i \text{ in spec. (R)} \; / \; |R/P_i|_r = |R|_r \}$ and

$\Lambda'(R) = \{ Q_i \text{ in spec. (R)} \; / \; |R/Q_i|_l = |R|_l \}$ .

If there is no confusion regarding the underlying ring R we may write $\Lambda$ and $\Lambda'$ for $\Lambda(R)$ and $\Lambda'(R)$ respectively .

(ii) For an ideal I of R we write ;

$\Lambda(R/I) = \{ (p/I) \text{ in spec. (R/I) } / \; p \; \varepsilon \; \text{spec.}(R) \text{ with } I \leq p \text{ and } |R/p|_r = |R/I|_r \}$ and

$\Lambda'(R/I) = \{ q/I \text{ in spec. (R/I)} \; / \; q \; \varepsilon \; \text{spec.}(R) \text{ with } I \leq q \text{ and } |R/q|_l = |R/I|_l \}$ .

**Remark(1) :-** It is not difficult to see that if $p/I \; \varepsilon \; \Lambda(R/I)$, then p/I is a minimal prime ideal of the ring R/I and hence P/I $\varepsilon$ min. spec.(R/I) . The same holds true for members of $\Lambda'(R/I)$ .

**Definition (2.2) :-** Let R be a Noetherian ring with nilradical N and let $\Lambda$ , $\Lambda'$ , A , B be as defined in definitions (2.1) and (1.1) respectively . Set $P = \cap \; p_i$ ; $p_i \; \varepsilon \; \Lambda$ and $Q = \cap \; q_i$ ; $q_i \; \varepsilon \; \Lambda'$ . Then the ideals P and Q are called the right and the left krull homogenous parts respectively of the nilradical of the ring R . Note that the factor rings R/P and R/Q are semiprime right and left krull homogenous rings respectively .

We now introduce the definition of a right ( or a left ) weakly krull symmetric ring , a weakly krull symmetric and a weakly ideal krull symmetric noetherian ring R .

**Definition (2.3) :-** Let R be a Noetherian ring with nilradical N and let $\Lambda$, $\Lambda'$, A, B be as defined in definitions (2.1) and (1.1) respectively. Let $P = \cap p_i$ ; $p_i \varepsilon \Lambda$ and let $Q = \cap q_i$ ; $q_i \varepsilon \Lambda'$. Then we define the following terms,

(i) R is called a right (or a left) weakly krull symmetric ring if $\Lambda'(R) \leq \Lambda(R)$ (or $\Lambda(R) \leq \Lambda'(R)$).

(i) R is called a weakly krull symmetric ring if $\Lambda(R) = \Lambda'(R)$.

(ii) R is called a weakly ideal krull symmetric ring if A = B.

**Proposition (2.4) :-** Let R be a Noetherian ring with nilradical N and let $\Lambda$, $\Lambda'$, A, B be as defined in definitions (2.1) and (1.1) respectively. Let $P = \cap p_i$ ; $p_i \varepsilon \Lambda$ and let $Q = \cap q_i$ ; $q_i \varepsilon \Lambda'$. Then the following conditions are equivalent ;

(i) R is a right weakly krull symmetric ring (that is $\Lambda'(R) \leq \Lambda(R)$).

(ii) For any ideal J of R $|R/J|_r < |R|_r$ implies that $|R/J|_l < |R|_l$.

(iii) For any ideal I of R $|R/I|_l = |R|_l$ implies that $|R/I|_r = |R|_r$.

We leave the statement of the analogous proposition for a left weakly krull symmetric ring to the reader. For weakly krull symmetric noetherian rings we have the following proposition ;

**Proposition (2.5) :-** Let R be a Noetherian ring with nilradical N and let $\Lambda$, $\Lambda'$, A, B be as defined in definitions (2.1) and (1.1) respectively. Let $P = \cap p_i$; $p_i \in \Lambda$ and let $Q = \cap q_i$; $q_i \in \Lambda'$. Then the following conditions are equivalent;

(i) R is a weakly krull symmetric ring (that is $\Lambda'(R) = \Lambda(R)$).

(ii) For any ideal J of R $|R/J|_r < |R|_r$ if and only if $|R/J|_l < |R|_l$

(iii) For any ideal I of R $|R/I|_l = |R|_l$ if and only if $|R/I|_r = |R|_r$.

**Proposition (2.6) :-** Let R be a Noetherian ring with nilradical N and let $\Lambda$, $\Lambda'$, A, B be as defined in definitions (2.1) and (1.1) respectively. Let $P = \cap p_i$; $p_i \in \Lambda$ and let $Q = \cap q_i$; $q_i \in \Lambda'$. Then the following conditions are equivalent;

(i) R a weakly ideal krull symmetric ring (that is A = B).

(ii) For any ideal J of R $|J|_r < |R|_r$ if and only if $|J|_l < |R|_l$.

(iii) For any ideal I of R $|I|_r = |R|_r$ if and only if $|I|_l = |R|_l$.

We state the following Remark in respect of our definitions;

**Remark (2) :-** No example is known of a noetherian ring R for which A ≠ B or for which $\Lambda(R) \neq \Lambda'(R)$. We may note that in case $\Lambda(R) = \Lambda'(R)$, then we have that P = Q. The main theorem of [2] guarantees that $\Lambda(R) = \Lambda'(R)$ in case R is both a left and a right krull homogenous ring. We state this fact in the theorem given below.

**Theorem (2.7) :-** Let R be a Noetherian ring that is a left and a right k-homogeneous ring (that is A = B = 0). Then we have that $\Lambda(R) = \Lambda'(R) = N(R)$, where $N(R)$ is the nilradical of R.

We now present the following proposition below ;

**Proposition (2.8)** :- Let R be a Noetherian ring with nilradical N and let $\Lambda$, $\Lambda'$, A, B be as defined in definitions (2.1) and (1.1) respectively. Let $P = \cap p_i$ ; $p_i \in \Lambda$ and let $Q = \cap q_i$ ; $q_i \in \Lambda'$. Then the following hold true ;

(a) $\Lambda(R) = \Lambda(R/A)$.

(b) $\Lambda'(R) = \Lambda'(R/B)$.

**Proof** :- (a) First observe that since $|R/A|_r = |R|_r$, hence, $\Lambda(R/A) \leq \Lambda(R)$. Now, let $P \in \Lambda(R)$. We claim that $P \in \Lambda(R/A)$.

Proof of the claim :- If $A = 0$, then there is nothing to prove. So, assume $A \neq 0$. Then using lemma (1.2) it follows that $r(A) \neq 0$ and $|R/r(A)|_r < |R|$. Hence, $r(A) \not\subset P$. Since, $Ar(A) = 0$ implies $Ar(A) \subseteq P$ so we get that $A \subseteq P$. Since $|R/A|_r = |R|_r$, it thus follows that $P \in \Lambda(R)$ implies that $P \in \Lambda(R/A)$. This proves the claim. Thus, $\Lambda(R) \leq \Lambda(R/A)$.

Now we conclude the proof by observing that since $|R/A|_r = |R|_r$, it is then obvious that $\Lambda(R/A) \leq \Lambda(R)$, so $\Lambda(R) = \Lambda(R/A)$.

(b) The proof of (b) follows as a left analogue of the proof of (a) above. We thus leave its proof to the reader.

Using the above proposition we state and prove a mild generalization of the main theorem of [2].

**Theorem (2.9)** :- Let R be a Noetherian ring with nilradical N. Let P and Q be the right and left krull homogenous parts of N respectively. Then if R is a weakly ideal krull symmetric noetherian ring then R is a weakly krull symmetric ring. Moreover in this case then $P = Q$.

**Proof :-** Let $\Lambda$, $\Lambda'$, A, B be as defined in definitions (2.1) and (1.1) respectively. From definition (2.2) above we have that $P = \cap p_i$; $p_i \in \Lambda$ and $Q = \cap q_i$; $q_i \in \Lambda'$. Moreover given R is weakly ideal krull symmetric ring we get that $A = B$. Then to prove the theorem we consider the ring $R/A$. Since $A = B$, by proposition (1.3) above we get that $R/A$ is a noetherian right and a left krull homogenous ring and hence applying theorem (2.7) above we get that $\Lambda(R/A) = \Lambda'(R/A)$. Applying proposition (2.8) above we get that $\Lambda(R) = \Lambda'(R)$. Hence R is a weakly krull symmetric ring.

**Remark (3) :-** It is unclear whether or not the converse of theorem (2.8) above is true or not.

We now prove our main theorem below ;

**Theorem (2.10) :-** Let R be a Noetherian ring with nilradical N and let $\Lambda$, $\Lambda'$, A, B be as defined in definitions (2.1) and (1.1) respectively. Let $P = \cap p_i$; $p_i \in \Lambda$ and let $Q = \cap q_i$; $q_i \in \Lambda'$. Then the following hold true ;

(i) If R is a right krull homogenous ring then R is a right weakly krull symmetric ring ( that is $\Lambda'(R) \leq \Lambda(R)$).

**(ii)** If R is a left krull homogenous ring then R is a left weakly krull symmetric ring ( that is $\Lambda(R) \leq \Lambda'(R)$).

**Proof :- (i**) To prove the result we first assume that R is a right krull homogenous ring. Let now B be an ideal of R that has the usual meaning. Then using proposition (1.5) above we get that R/B is a left and a right krull homogenous ring with $|R/B|_r = |R|_r$ and $|R/B|_l = |R|_l$. Hence by theorem (2.7) we have that $\Lambda'(R/B) = \Lambda(R/B)$. By theorem (2.8) we have that $\Lambda'(R/B) = \Lambda(R/B)$. Thus we get that $\Lambda'(R) = \Lambda'(R/B) = \Lambda(R/B)$. Since $|R/B|_r = |R|_r$, so $\Lambda(R/B) \leq \Lambda(R)$. Hence we must have that $\Lambda'(R) \leq \Lambda(R)$.

**(ii)** The proof of (ii) follows on the same lines as that of (i) above.

We end this paper with a remark and an acknowledgement;

**Remark and acknowledgement :-** We remark that if R is a Noetherian ring, then R is called a krull symmetric ring if $|R|_r = |R|_l$. We must mention that the present paper is a replacement of two earlier papers of the author with the title "On the krull symmetry of a noetherian ring" that appeared in the Arxiv. In these papers the author claimed to have proved that a noetherian ring R is krull symmetric. But in several correspondences with Prof. T.H.Lenagan the author found these proofs to be incomplete. The author is much indebted to Prof. T.H.Lenagan for his comments and advice.

**References:-**

**(1) K.R. Goodaarl and R.B. Warfield , "An Introduction to Non-commutative Noetherain Rings" , L.M.S., student texts ,16.**

**(2) C.L. Wangneo ; On A Certain Krull Symmetry Of a Noetherian Ring. arXiv:1110.3175v2[math.RA],Dec.2011.**